\documentclass[11pt]{article}

\usepackage[centertags]{amsmath}
\usepackage{amsfonts}
\usepackage{amssymb}
\usepackage{theorem}
\usepackage{ulem}
\usepackage{epsfig}
\usepackage[all]{xy}
\usepackage{wasysym}
\usepackage{calc}
\usepackage{color}
\theoremstyle{break} \newtheorem{theorem}{Theorem}[section]
\theoremstyle{break} 
\theoremstyle{break} \newtheorem{definition}[theorem]{Definition}       
\theoremstyle{break} \newtheorem{lemma}[theorem]{Lemma}
\theoremstyle{break} \newtheorem{corollary}[theorem]{Corollary}
\theoremstyle{break} \newtheorem{example}[theorem]{Example}
\theoremstyle{break} 
\theoremstyle{break}
{\theorembodyfont{\rmfamily}\newtheorem{remark}[theorem]{Remark}}
{\theorembodyfont{\rmfamily}}
\theoremstyle{break} 
\theoremstyle{break} 
\theoremstyle{break} 
\theoremstyle{break} 
\numberwithin{equation}{section}

\newcommand{\R}{{\mathbb{R}}}
\newcommand{\D}{{\mathbb{D}}}
\newcommand{\C}{{\mathbb{C}}}

\def\Re{\mathop{{\rm Re}}}
\def\Im{\mathop{{\rm Im}}}

\newcommand{\E}{{\rm EU}}
\newcommand{\RI}{{\rm RI}}

\renewcommand{\phi}{\varphi}
\newcommand{\eps}{\varepsilon}

\begin{document}

\renewcommand{\thefootnote}{}
\stepcounter{footnote}
\begin{center}
{\bf \LARGE Beurling's free boundary value problem\\[2mm] in conformal geometry}\footnote{D.K.~was supported by a HWP scholarship. O.R.~received partial support from
the German--Israeli Foundation (grant G-809-234.6/2003). 
\hfill{{\it Israel Journal Math.} to appear}}
\end{center}
\renewcommand{\thefootnote}{\arabic{footnote}}
\setcounter{footnote}{0}
\begin{center}
{\large  \sc Florian Bauer, Daniela Kraus, Oliver Roth and Elias Wegert}\\[2mm]
\end{center}

\centerline{January 10, 2009}

\renewcommand{\thefootnote}{\arabic{footnote}}

\smallskip
\begin{center}
\begin{minipage}{13cm} {\bf Abstract.}   
The subject of this paper is Beurling's  celebrated extension 
of the Riemann mapping theorem \cite{Beu53}. Our  point of departure is the observation that  the only known
 proof of the Beurling--Riemann mapping theorem  contains 
a number of gaps which seem 
 inherent in Beurling's geometric and approximative approach.  We provide
 a complete proof of the Beurling--Riemann mapping theorem by combining Beurling's geometric method with a number of new analytic tools, notably $H^p$--space techniques and methods from the theory of 
Riemann--Hilbert--Poincar\'e problems. One additional advantage  of this approach is that it
leads to an extension of the Beurling--Riemann mapping theorem 
 for analytic maps  with prescribed branching. Moreover, 
it allows a complete 
description of the boundary regularity of solutions in the (generalized) Beurling--Riemann
 mapping theorem extending  earlier results that have been obtained by PDE techniques.   We finally consider the question of uniqueness in the extended Beurling--Riemann mapping theorem. 
\end{minipage}
\end{center}

\section{Introduction}

Let  $\mathcal{H}_0(\D)$ denote the set of 
analytic functions $f$ on the unit disk $\D:=\{z \in \C \, : \, |z|<1\}$
 normalized by $f(0)=0<f'(0)$ and let $\Phi : \C \to \R$ be a continuous, 
positive and bounded function.
 Beurling's conformal geometric free boundary value problem  \cite{Beu53} asks for
{\it univalent} functions $f \in {\cal H}_0(\D)$  that satisfy
\begin{equation} \label{eq:boundarycondition}
 \lim \limits_{z\to \xi } \left( |f'(z)|-\Phi(f(z)) \right)=0 \,,  \qquad \xi
 \in \partial \D .
\end{equation}
We call any $f \in \mathcal{H}_0(\D)$ (univalent or not) for which
(\ref{eq:boundarycondition}) holds a solution for $\Phi$.

\smallskip

Beurling's paper \cite{Beu53} and its successor \cite{Beu58}
proved to be quite influential in various different branches of mathematics
such as partial differential equations, geometric
function theory and Riemann--Hilbert problems.
For instance, some of Beurling's ideas
 are nowadays extensively used in the theory of
free boundary value problems for PDEs and in fact the papers \cite{Beu53,Beu58} are
widely  considered as some of the pioneering papers on free boundary value problems
(see  \cite{GS,H}). They also found considerable attention in geometric
complex analysis and conformal geometry, see for instance
\cite{AB96,Av94,Av96,FRold,FR01,HMT03,Ro07} as some of the more recent references.

\smallskip

One of the purposes of the present paper is to provide a thorough and complete 
discussion of Beurling's original free boundary value problem (\ref{eq:boundarycondition}). This requires  advanced analytic tools from
$H^p$--space theory and Riemann--Hilbert--Poincar\'e problems, which have not been used for this purpose before. These tools make it possible to 
 establish in addition  a number of extensions of 
Beurling's results e.g.~to solutions  with prescribed critical points,  to describe the
boundary behaviour of the solutions and to deduce new sufficient conditions for uniqueness of solutions.

\smallskip

Beurling's treatment of the boundary value problem
(\ref{eq:boundarycondition}) in \cite{Beu53}
 is based on an ingenious geometric version of Perron's
method from the theory of subharmonic functions. He defines a class of
supersolutions and subsolutions for $\Phi$ (see below for the precise
definition) and shows that there is always 
a ``largest'' univalent subsolution $f^*$  and a ``smallest'' univalent
supersolution $g^*$. He then asserts that both $f^*$ and $g^*$ are
in fact solutions for $\Phi$, but his arguments are in both cases incomplete.
For instance,  Beurling's approach  to prove that
 $g^*$  is  a
  solution for $\Phi$ is in two steps. First he deals with the special case that  $g^*(\D)$ is of  ``Schoenfliess type''. In order to handle the general case,
 he then makes use of the assertion that every strictly shrinking sequence of  (normalized) 
  simply connected domains converges in the sense of kernel convergence to a
  domain of Schoenfliess type. This, however, is not
  true in general (see Appendix 2 below for an explicit counterexample), so Beurling's method is destined to fail here. 
In order to circumvent this difficulty we combine  Beurling's geometric
 approach with more advanced analytic tools.
 As a result, we obtain an efficient  
  method  that avoids  ``domains of Schoenfliess type'' altogether and allows 
   a  treatment of the smallest univalent supersolution without
   unnecessary approximation techniques.

\smallskip

We also establish the existence of {\it nonunivalent}
  solutions for $\Phi$ with {\it prescribed branch points}. There are a number of
  reasons for taking nonunivalent solutions into account. For instance, 
  it seems  indispensable to find  first nonunivalent solutions
  in order to prove that there are always univalent solutions.
   The proof of existence of nonunivalent solutions in this paper
 is based on a fixed--point argument and closes another possible gap in 
 Beurling's original approach, see Section \ref{par:gen}.
  A second reason for allowing nonunivalent solutions  is that 
  Beurling's problem might be viewed as a special case of a certain type of generalized
  Riemann--Hilbert problems (sometimes called Riemann--Hilbert--Poincar\'e
  problems). These problems deal with the construction of analytic maps with
  prescribed boundary behaviour {\it and} preassigned branch points.
   We note that in this context a variant of Beurling's boundary value problem with specified critical
  points played a key role in recent work on ``hyperbolic'' Blaschke products
  (see \cite{FRold}) and on infinite Blaschke products with
  infinitely many critical points (see \cite{KR}).

\smallskip

Finally, we  discuss in detail the boundary behaviour
of the solutions for $\Phi$. 
For {\it univalent} solutions  Beurling's problem  (\ref{eq:boundarycondition}) 
might be viewed as  a free boundary value problem for  $\Omega=f(\mathbb{D})$
involving PDEs (see Appendix 1), but this relation breaks down for
nonunivalent solutions.
Thus one can apply techniques from PDEs to
study the boundary regularity of univalent solutions, but not for nonunivalent solutions.
In particular, it follows from  results of 
Alt \& Caffarelli \cite{AC}, Kinderlehrer \& Nirenberg
\cite{KN77} and Gustafsson \& Shahgholian \cite{GS} that every
univalent solution $f$ for $\Phi$ is of class
$C^{1,\alpha}(\partial \D)$ for {\it some} $0<\alpha<1$ when $\Phi$ is 
H\"older continuous, that $f \in C^{k+1,\alpha}(\partial \D)$ when 
$\Phi \in C^{k,\alpha}$, $k \ge 1$ and $0<\alpha<1$, and that
$f$ is real analytic across $\partial \D$ when $\Phi$ is real analytic (see
also Sakai \cite{S1,S2}). The result for real analytic $\Phi$
has recently been extended to all (i.e.~not necessarily
univalent) solutions in \cite{Ro07} using a completely different approach. 
Based on a method specific to Beurling's problem,
we complement the results of   \cite{AC,KN77,GS}
in Theorem \ref{thm:regularity} below  by showing that
every solution for $\Phi$ belongs
to $C^{k+1,\alpha}(\partial \D)$ provided $\Phi$ is of class $C^{k,\alpha}$
for
 $0<\alpha<1$ and {\it all} $k \ge 0$.

\smallskip

This paper is organized as follows. We start in Section \ref{sec:subsolutions}
with a discussion of the set $A_{\Phi}$ 
of subsolutions to Beurling's boundary value problem
(\ref{eq:boundarycondition}) including a number of simplifications and generalizations of Beurling's original treatment of $A_{\Phi}$.  In Section \ref{sec:supersolutions} we show that the set $B_{\Phi}$ of {\it univalent} supersolutions can be handled almost
identically as the set $A_{\Phi}$ using the (usual) topology of locally
uniform convergence on the unit disk. In order to incorporate the boundary
behaviour of univalent supersolutions we continuously embed $B_{\Phi}$ in the
Hardy space $H^p$, $0<p<1/2$. This is 
possible by results of Feng and
MacGregor \cite{FM76} on the integral means of the derivative of univalent
functions and Hardy--Littlewood--type arguments. In this context, the key result
is  Lemma \ref{lem:H^p} below.
Section \ref{sec:solutions} is divided into four parts.
In \S \ref{par:gen} we consider a class of Riemann--Hilbert--Poincar\'e
problems, which includes Beurling's free boundary value
problem as special case. 
We first establish a representation formula for  solutions to this
more general type of boundary value problem. The existence of such solutions is then proved
by an application of Schauder's fixed point principle.  
 The representation formula  also leads  to
a full description of the boundary behaviour of solutions in \S \ref{par:reg}.
Armed with at least nonunivalent solutions we
 then return to Beurling's original boundary problem in \S
\ref{par:maxi} and \S\ref{par:mini}. In \S
\ref{par:maxi} we will see that the maximal subsolution  is in fact
a (univalent) solution. Here we make essential use of the results of \S \ref{par:gen}. 
The proof that the minimal univalent  supersolution  is 
a solution is much more elaborate and is given in \S
\ref{par:mini}. In  Section
\ref{sec:uniqueness} we briefly discuss the question of uniqueness in
Beurling's boundary value problem and find (slight) generalizations
of uniqueness results due to Beurling \cite{Beu53} and Gustafsson \& Shahgholian \cite{GS} (see also \cite{HMT03}).
We conclude the paper with two appendices. Appendix 1 indicates how Beurling's
 problem for univalent functions is connected with a class of free
boundary value problems for PDEs, which also arise in many areas of physics 
(Hele--Shaw flows) as well as in mathematical analysis (Quadrature
domains). Finally, in Appendix 2 we discuss an explicit counterexample to
Beurling's method of proof in \cite{Beu53}.

\section{Subsolutions} \label{sec:subsolutions}

In the sequel we
denote the Poisson kernel on $\D$ by
$$ P(z,e^{i t}):=\Re \left(\frac{e^{it}+z}{e^{it}-z} \right) \, , \qquad z \in
\D\,,\, \,  t  \in \R \, ,$$
and use the notation
$$ ||\varphi||:=\sup \limits_{w \in \C} |\varphi(w)|$$
for any bounded function $\varphi : \C \to \C$.

\begin{definition}
Let $\Phi:\C \to \R$ be a positive, continuous and bounded function. The
set $A_{\Phi}$ is defined by
$$ A_{\Phi}:= \left\{ f \in {\cal H}_0 (\D) \, : \, \limsup_{|z| \to 1} \left(
   |f'(z)|- \Phi(f(z)) \right) \le 0\right\}\, .$$
We call every function $f \in A_{\Phi}$ a subsolution for $\Phi$.
\end{definition}

The goal of this section is to show that there is always a unique ``largest'' 
subsolution $f^*$ for $\Phi$ and that this largest subsolution is {\it
  univalent}:

\begin{theorem}\label{thm:maxi_function}
Let $\Phi:\C \to \R$ be a positive, continuous and bounded function.
Then there exists a unique  function $f^{*} \in A_{\Phi}$ such that 
$$f^{* \, \prime}  (0)= \sup \limits_{f \in A_{\Phi}} f'(0) \, .$$
The function $f^*$ is  univalent in $\D$ with $f^*( \D)=A^*$, where
$$ A^*:=\bigcup  \limits_{ f \in A_{\Phi}} f(\D)\, .$$
In particular, $f(\D)\subseteq A^*$ for all $f \in A_{\Phi}$.
\end{theorem}

We call the function $f^*$ of Theorem \ref{thm:maxi_function} the maximal
univalent subsolution for $\Phi$.
The following simple facts about subsolutions will be used in the proof of
Theorem \ref{thm:maxi_function}.

\begin{lemma} \label{lem:extended0}
Let $\Phi : \C \to \R$ be a positive, continuous and bounded function.\\[-6mm]
\begin{itemize}
\item[(a)] Any subsolution for $\Phi$ has a (Lipschitz) continuous extension to the closed
  unit disk $\overline{\D}$. The set $A_{\Phi}$ is uniformly bounded on
  $\overline{\D}$ and equicontinuous at every point $z_0 \in \overline{\D}$.

\item[(b)] A function $f \in \mathcal{H}_0(\D)$ with a continuous extension to
  $\overline{\D}$ is a subsolution for $\Phi$
  if and only if 
\begin{equation} \label{eq:Aphi1}
 \log |f'(z)| \le \frac{1}{2 \pi} \int \limits_{0}^{2 \pi} P(z,e^{i t}) \, \log
\Phi(f(e^{it})) \, dt\, , \qquad z \in \D \, .  
\end{equation}
\item[(c)] If a sequence of subsolutions for $\Phi$ converges locally
  uniformly in $\D$, then the limit function $f$ is either
  again a  subsolution for $\Phi$ or $f \equiv 0$.
\item[(d)] Let $f \in A_{\Phi}$ and let
$\Psi : \C \to \R$ be a positive, continuous and bounded
  function with $\Psi > \Phi$.
  Then for all $r<1$ sufficiently close to $1$  the
  function $f_r(z):=f(rz)$ is a subsolution for $A_{\Psi}$.
\end{itemize}
\end{lemma}

\begin{remark}
 Lemma \ref{lem:extended0} shows that for treating the class $A_{\Phi}$
 one can use the uniform topology on the {\it
   closed} unit disk $\overline{\D}$. This facilitates the handling of
 subsolutions.  
Beurling proved Lemma \ref{lem:extended0} (b) and (c) for {\it univalent}
subsolutions. We need the general case. \end{remark}

{\bf Proof.} 
(a)\,  The maximum principle implies
that for any $f \in A_{\Phi}$ the function $|f'|$ 
is bounded above in $\D$ by $M:=||\Phi||$, so $f$ has 
a (Lipschitz) continuous extension to $\overline{\D}$. Note that the  Lipschitz
constant is independent of $f$, so $A_{\Phi}$ is uniformly bounded and equicontinuous in $\overline{\D}$.

\smallskip

(b) \,  This follows from the fact that
the Poisson integral on the righthand side of (\ref{eq:Aphi1}) is
  harmonic in $\D$ with boundary values $\log \Phi(f(\xi))$, $|\xi|=1$.

\smallskip

(c) \, Let $f_n \in A_{\Phi}$ and suppose that $(f_n)_n$ converges locally
uniformly to a holomorphic function $f : \D \to \C$ with $f(0)=0$. 
Then by (a) and the Arzel\`a--Ascoli theorem a subsequence
$(f_{n_k})_k$ converges uniformly on $\overline{\D}$, so $f : \overline{\D}
\to \C$ is continuous. If $f \not \equiv 0$, then $f$ is not constant. Hence
$\log |f'(z)|$ is a subharmonic function on $\D$ and part (b) shows that
 $f \in A_{\Phi}$ because
\begin{eqnarray*} \log |f^{\prime}(z)|& =& \lim \limits_{k \to \infty} \log 
|f_{n_k}'(z)| \le
\lim \limits_{k \to \infty} \frac{1}{2 \pi} \int
  \limits_{0}^{2 \pi} P(z,e^{it}) \, 
\log \Phi(f_{n_k}(e^{it}))
 \, dt \\ &=& \frac{1}{2 \pi} \int
  \limits_{0}^{2 \pi} P(z,e^{it}) \, 
\log \Phi(f(e^{it})) \, dt \, , \qquad z \in \D \, . 
\end{eqnarray*}

(d) \,  Assume to the contrary that $f_{r_n} \not \in A_{\Psi}$ for some
  sequence of radii $r_n \nearrow 1$. Then there are points $|\xi_n|=1$ such that
$r_n |f'(r_n \xi_n)|=
|f^{\prime}_{r_n}(\xi_n)| > \Psi(f_{r_n}(\xi_n))=\Psi(f(r_n \xi_n)) \, . $
We may assume $r_n \xi_n \to \xi \in \partial \D$. Then
$$ \limsup_{z \to \xi} |f'(z)| \ge \limsup_{n \to \infty} |f'(r_n \xi_n)| \ge
\Psi(f(\xi))>\Phi(f(\xi)) \, , 
$$
which contradicts $f \in A_{\Phi}$. \hfill{$\blacksquare$}

\medskip

The main tool needed for the proof of Theorem \ref{thm:maxi_function} is a
geometric version of the Poisson modification of a subharmonic function.
For this Beurling introduced an ``extended'' union of domains in the complex
plane which is always  simply connected.

\begin{definition}[Extended union]
Let $D_1$ and $D_2$ be two bounded domains in $\C$ with $0 \in D_1 \cap D_2$. 
We call the
complement of the unbounded component of $\hat{\C} \backslash (D_1 \cup D_2)$ the
extended union of $D_1$ and $D_2$ and denote it by $\E(D_1,D_2)$. 
\end{definition} 

Beurling has a formally different, but equivalent definition for $\E(D_1,D_2)$.
We note the following easily verified properties of the extended union.\\[-5.5mm]
\begin{itemize}
\item[]{\bf(EU1)}  \, $\E(D_1,D_2)$ is the smallest simply connected domain which contains $D_1
\cup D_2$.
\item[]{\bf (EU2)} \, If $D_1 \subseteq D_1'$ and $D_2 \subseteq D_2'$, then
$\E(D_1,D_2) \subseteq \E(D_1',D_2')$.
\item[]{\bf (EU3)} \, $\partial \E(D_1,D_2) \subseteq \partial D_1 \cup \partial D_2$. 
\end{itemize}

\begin{definition}[Upper Beurling modification]
Let $f_1, f_2 \in \mathcal{H}_0(\D)$. Then the  conformal
map $f$ from $\D$ onto $\E(f_1(\D),f_2(\D))$ normalized by $f(0)=0$ and
$f'(0)>0$ is called the upper Beurling modification of $f_1$ and $f_2$.
\end{definition}

\begin{lemma}\label{lem:upper_mod}
Let $\Phi : \C \to \R$ be a positive, continuous and bounded function and
$f_1, f_2$  two subsolutions for $\Phi$. Then the upper Beurling
modification of $f_1$ and $f_2$ is also a subsolution for $\Phi$.
\end{lemma}

The following proof is different from Beurling's proof insofar as we replace
Beurling's Riemann surface construction by a simple application of the
Julia--Wolff lemma (see \cite[p.~82]{Pom92}).

\medskip

{\bf Proof.}
(i) \, 
We first prove the lemma under the additional assumption that $f_1$ and $f_2$ are
analytic on $\overline{\D}$. In this case the boundary of the extended union
$D:=\E(f_1(\D),f_2(\D))$ is locally connected and
consists of finitely many analytic Jordan arcs. Thus the upper Beurling
modification $f$ of $f_1$ and $f_2$
has an analytic continuation across
$\partial \D \backslash N$, where $N \subset \partial \D$ is  a finite set, and
  $D$ is a Smirnov domain
\cite[p.~60 \& p.~163]{Pom92}, i.e., $f$ satisfies 
\begin{equation} \label{eq:extended1}
 \log |f'(z)|=\frac{1}{2\pi} \int \limits_{0}^{2 \pi} P(z,e^{it}) \, \log
|f'(e^{ it})| \, dt \, , \qquad z \in \D \, . 
\end{equation}
Since $\partial D \subset f_1(\partial \D) \cup f_2(\partial \D)$,  there is
in particular for every $\xi \in \partial \D\backslash N$ a point $\xi_0 \in \partial \D$
such that $f(\xi)=f_1(\xi_0)$ or $f(\xi)=f_2(\xi_0)$.
If $f(\xi)=f_1(\xi_0)$, then the 
function $w(z):=f^{-1}(f_1(z))$ maps $\D$ into $\D$ with $w(0)=0$ and
$w(\xi_0)=\xi$. By the Julia--Wolff lemma $|w'(\xi_0)| \ge 1$, so
\begin{equation} \label{eq:extended2}
|f'(\xi)| =\frac{|f_1'(\xi_0)|}{|w'(\xi_0)|} \le |f_1'(\xi_0)| \le
\Phi(f_1(\xi_0)) =\Phi(f(\xi)) \, .
\end{equation}
The same conclusion holds if $f(\xi)=f_2(\xi_0)$.
Thus (\ref{eq:extended2}) holds for every $\xi \in \partial \D$ except for
finitely many points. Therefore (\ref{eq:extended1}) leads to
$$ \log |f'(z)|  \le  \frac{1}{2\pi} \int \limits_{0}^{2 \pi} P(z,e^{it}) \, \log
\Phi(f(e^{ it})) \, dt \, , \qquad z \in \D \, . $$
This shows $f \in A_{\Phi}$.

\smallskip

(ii) \, We now turn to the general case. Let $\Psi : \C \to \R$ be a positive, continuous and bounded
  function with $\Psi > \Phi$. In view of Lemma \ref{lem:extended0} (d),  the functions $f_1( rz)$ and $f_2(rz)$ belong to
$A_{\Psi}$ for all $0<r_0<r<1$. By what we have shown above, the upper
Beurling modification $f_r$ of $f_1(rz)$ and $f_2(rz)$ belongs to $A_{\Psi}$
for every $0<r_0<r<1$. Now $f_r$ maps $\D$ conformally onto
$D_r:=\E(f_1(r \D),f_2(r \D))$. Note that $D_r \subseteq D_{r'}$ whenever $0<r<r'<1$ in
view of (EU2). Thus  
as $r \to 1$ the domains $D_r$ converge in the
sense of kernel convergence to the simply connected domain
$$D':=\bigcup \limits_{r_0<r<1} D_r\,  $$
with $f_1(\D) \cup f_2(\D) \subset D'$. This implies $D \subset D'$ by (EU1).
On the other hand, by (EU2), we have $D_r \subseteq D$, so $D' \subseteq D$.
Carath\'eodory's convergence theorem shows  $f_r \to f$ locally uniformly in
$\D$, see \cite[Chap.~1.4]{Pom92}. Since $f_r \in A_{\Psi}$ for all $r$ close enough to $1$, 
Lemma \ref{lem:extended0} (c) gives $f \in A_{\Psi}$. As $\Psi$ is an
arbitrary positive, continuous and bounded
  function on $\D$ with $\Psi > \Phi$, we conclude $f \in A_{\Phi}$.
\hfill{$\blacksquare$}

\begin{remark}
The approximation argument in the above proof cannot be avoided entirely.
This is due to the fact that subsolutions even though they are Lipschitz continuous up
to the unit circle may have a very bad behaved derivative. For instance 
the well--known example of Duren, Shapiro and Shields
\cite{DSS} (see also \cite[p.~159]{Pom92}) of a conformal map $f$ not of Smirnov
type satisfies $|f'|<1$ in $\D$, so belongs to $A_{\Phi}$ for
$\Phi \equiv 1$. Even for a solution $f$ for $\Phi$ the derivative $f'$ does not
need to have a continuous extension to $\overline{\D}$,
 see Example \ref{ex:pom}.
\end{remark}

{\bf Proof of Theorem \ref{thm:maxi_function}.} 
By Lemma \ref{lem:extended0} there is a function 
$f^* \in A_{\Phi}$ such that
$$ f^{*\, \prime}(0)=\max \limits_{f \in A_{\Phi}} f'(0) \, .$$
We need to prove that $f^*$ is univalent and $f(\D)\subseteq f^*(\D)$ for every $f \in
A_{\Phi}$. Let $f \in A_{\Phi}$ and let $F$ be the upper Beurling modification
of $f$ and $f^*$.  Lemma \ref{lem:upper_mod} shows $F \in A_{\Phi}$, so
$F'(0) \le f^{* \, \prime}(0)$. On the other hand, we have 
$F(\D)=\E(f^*(\D),f(\D)) \supseteq f^*(\D)$. By the principle of subordination
this implies $F=f^*$. Hence $f^*$ is univalent and $f(\D) \subseteq f^*(\D)$. 
Now it is also clear that $f^*$ is uniquely determined.
\hfill{$\blacksquare$}

\section{Univalent supersolutions} \label{sec:supersolutions}

\begin{definition}
Let $\Phi:\C \to \R$ be a positive, continuous and bounded function. Then the
set $B_{\Phi}$ is defined by
$$ B_{\Phi}:= \left\{ g \in {\cal H}_0 (\D) \, : \, g \text{ univalent and } \liminf_{|z| \to 1} \left(
   |g'(z)|- \Phi(g(z)) \right) \ge 0\right\}\, .$$
We call any function $g \in B_{\Phi}$ a univalent supersolution for $\Phi$.
\end{definition}

Note that we  consider only {\it univalent} supersolutions, see Remark
\ref{rem:univalent_supersol} for a partial explanation.
We establish in this paragraph the following counterpart
to Theorem \ref{thm:maxi_function}.

\begin{theorem}\label{thm:mini_function}
Let $\Phi:\C \to \R$ be a positive, continuous and bounded function.
Then there
exists a unique   function $g^{*} \in B_{\Phi}$ such that
$$g^{* \, \prime}  (0)= \inf \limits_{g \in B_{\Phi}} g'(0) \, . $$
The function $g^*$ maps $\D$ conformally onto $B^*$, where
$$B^*:=\bigcap  \limits_{ g \in B_{\Phi}} g(\D)\, .$$
In particular, $B^* \subseteq g(\D)$ for all $g \in B_{\Phi}$.
\end{theorem}

We call the function $g^*$ of Theorem \ref{thm:mini_function} the minimal
univalent supersolution for $\Phi$.

\smallskip

Supersolutions are considerably more complicated than subsolutions since
they do not need to be continuous up to
the unit circle. In particular, the topology of uniform convergence
on the closed unit disk is inappropriate for  the class $B_{\Phi}$.  
Nevertheless,  $B_{\Phi}$ can  
 be handled in a
similar way as the class $A_{\Phi}$, but now  using the notion of locally uniform
convergence in $\D$. For treating Beurling's {\it boundary
  value problem} we ultimately need to pass from inside the unit disk to the unit
circle.

\begin{lemma}\label{lem:B_Phi} 
Let $\Phi:\C \to \R$ be a positive, continuous and bounded function.
\begin{itemize}
\item[(a)] 
Any univalent supersolution $g$ for $\Phi$ has radial limits almost
  everywhere, i.e.~the limit 
$$g(e^{i t})=\lim \limits_{r \to 1} g(r e^{i  t}) \, , $$
exists for a.e.~$e^{ it} \in \partial \D$. 
  The boundary function $g$ belongs to $L^p(\partial \D)$ for
  every $0<p<1/2$. If $g$ is bounded, then $g^{-1}$ has a continuous extension
  to the closure of $g(\D)$.
\item[(b)]
A bounded and univalent function $g \in {\cal H}_0 (\D)$ belongs to $B_{\Phi}$
 if and only if
\begin{equation} \label{eq:ex1}
 \log|g'(z)|\ge \frac{1}{2\pi} \int_0^{2\pi} P(z, e^{it}) \log
\Phi(g(e^{it}))\, dt\, , \qquad z \in \D\, .
\end{equation}
\item[(c)]
If a uniformly bounded sequence  of univalent supersolutions for $\Phi$ converges locally uniformly
in $\D$, then the limit function $g$ is again a univalent supersolution for $\Phi$.
\item[(d)] Let $\Psi : \C \to \R$ be a positive, continuous and bounded
  function with $\Psi < \Phi$ and let  $g $ be a univalent supersolution for
$\Phi$. If $g$ is bounded, then  for all $r<1$ sufficiently close to $1$ the function
$g_r(z):=g(rz)$ is a univalent supersolution for $\Psi$. 
\end{itemize}
\end{lemma}

{\bf Proof of Lemma \ref{lem:B_Phi} (a) $\&$ (b) $\&$ (d):}

\smallskip

(a)\,
A univalent holomorphic function in $\D$ belongs to the Hardy spaces $H^p$
for $0<p<1/2$ and has therefore radial limits almost everywhere, see
\cite{Dur2000}. If $g$ is bounded, then 
$$|g'(z)| \ge c:=\inf_{|w| \le ||g||}
\Phi(w)>0 \, , $$  
so $|(g^{-1})'(w)| \le 1/c$ in $g(\D)$, i.e.~$g^{-1}$ has
a continuous extension to the closure of $g(\D)$.

\smallskip

(b)\,  First note that if $g \in {\cal H}_0(\D)$ is a bounded 
and univalent function, then, since $g(\D)$ is simply connected, there is a unique
solution $U_g$ to the Dirichlet problem
$$
\begin{array}{rclcr}
\Delta u& \equiv& 0 & \text{in} & g(\D) \, \phantom{.}\\[1mm]
u& = &\log \Phi & \text{on} & \partial g(\D)\, .
\end{array}
$$
The function $U_g \circ g$ is harmonic in $\D$ and  in view of
part (a) has the radial limit
$$(U_g\circ g)  (e^ {it})=\log \Phi (g(e^{it})) \quad  $$
for a.e.~$e^{ it} \in \partial \D$.
Hence we can write
\begin{equation}\label{eq:pois}
(U_g\circ g) (z)= \frac{1}{2\pi} \int_0^{2\pi} P(z, e^{it}) \log
\Phi(g(e^{it}))\, dt\, , \quad   z \in \D\, .
\end{equation}
Now let $g \in {\cal H}_0(\D)$ be univalent and bounded. Then $g \in B_{\Phi}$ 
 if and only if 
$$ \liminf_{|z| \to 1} \left( \log|g'(z)| - \log \Phi(g(z))   \right) \ge 0\,
.$$
By the univalence of $g$ we conclude  that the latter inequality is equivalent to
$$ \liminf_{|z| \to 1} \left( \log|g'(z)| - U_g(g(z))   \right) \ge 0\, .$$
The minimum principle of harmonic functions  implies now that
this is  the same as
$$  \log|g'(z)| - U_g(g(z))    \ge 0\, ,\quad z \in \D\, . $$
 Thus,  by equation (\ref{eq:pois}), $g \in B_{\Phi}$ if and only if (\ref{eq:ex1})
holds.

\smallskip

(d) The proof is similar to the proof of Lemma \ref{lem:extended0} (d) taking into account that $g$ is bounded and that we may assume that $g(r_n \zeta_n) \to \eta \in \partial g(\D)$. The details are therefore omitted.
\hfill{$\blacksquare$}

\medskip

Part (c)  of Lemma \ref{lem:B_Phi} is much more difficult to prove. Its statement is  essentially due to Beurling \cite{Beu53}, but he does not
provide a proof. 
We first show that
a locally uniformly convergent sequence of 
 univalent holomorphic functions converges in a weak sense
also on the boundary. 

\begin{lemma}\label{lem:H^p}
Let $(f_n)_n \subset {\cal{H}}_0(\D)$ be a sequence of univalent functions.
 Then for any $0<p<1/2$ and any $f \in H^p$ the following are equivalent:
\begin{itemize}
\item[(i)]
$(f_n)_n$ converges to $f$ locally uniformly in $\D$.
\item[(ii)] $\displaystyle \lim \limits_{n \to \infty}  \int \limits_0^{2\pi} |f_n(e^{i
    t}) -f( e^{it})|^p \, d t =0$.
\end{itemize}
\end{lemma}

Lemma \ref{lem:H^p} implies that the class ${\cal S}=\{ f \in {\cal H}_0(\D): f
  \text{ univalent and } f'(0)=1  \}$ is {\it compactly} contained in $H^p$ for any $0<p<1/2$.

\medskip

{\bf Proof.}
The statement (ii) $ \Rightarrow$ (i) follows directly from the estimate
$$|f(z)|\le 2^{1/p} \left( \frac{1}{2 \pi} \int_{0}^{2 \pi} |f( e^{it})|^p
\, d t \right)^{1/p} \, (1-|z|)^{-1/p}\, , \quad z \in \D,$$
which holds for any $f \in H^p$, $0<p<\infty$ ; see
\cite[p.~36]{Dur2000}.

\smallskip

To prove the implication (i) $ \Rightarrow$ (ii) we follow the proof  of Theorem 5 in \cite{Gwi36}. Fix $0< \rho < 1$ and choose $r$ such that $ \rho < r<1$. 
Then 
\begin{equation*}
\begin{split}
 & \int \limits_{0}^{2 \pi} |f_n(r e^{i t}) -f(r e^{i t})|^p \,
dt  \\[2mm]
& \qquad \le \int \limits_{0}^{2 \pi} |f_n(r e^{i t}) -f_n(\rho e^{i t})|^p \,
dt +
\int \limits_{0}^{2 \pi} |f_n(\rho e^{i t}) -f(\rho e^{i t})|^p \,
dt+
\int \limits_{0}^{2 \pi} |f(r e^{i t}) -f(\rho e^{i t})|^p \,
dt\, .
\end{split}
\end{equation*}

The first and  the third integral have the same structure and can therefore be
handled simultaneously. We may assume that $f$ is univalent, because otherwise
the third integral vanishes. 
Let $r_k=r(1-1/2^k)$ and let $N$ be the uniquely determined integer with $r_N \le \rho < r_{N+1}$. Denote by $\varphi$ one of the functions $f_n$, $n=1,2, \ldots$, 
or $f$.
Then, by Hardy--Littlewood (see \cite[Thm.~1.9]{Dur2000}), 
\begin{eqnarray*}
\int \limits_{0}^{2 \pi} |\varphi(r e^{i t}) -\varphi(\rho e^{i t})|^p \,
dt   & \le &
\int \limits_{0}^{2 \pi} \left(\, \,  \int \limits_{r_N}^r |\varphi'(x e^{i t})| \,
  dx \right)^p dt \le
\sum_{k=N}^\infty \int \limits_{0}^{2 \pi} \left( \int \limits_{r_k}^{r_{k+1}} |\varphi'(x e^{i t})| \,
  dx \right)^p dt  \\[2mm]
 & \le & \sum_{k=N}^\infty (r_{k+1} -r_k)^p \int_0^{2 \pi} \left( \max_{r_k \le x \le r_{k+1}}
 |\varphi'(x e^{i t})| \right)^p \, d t\\[2mm]
& \le & C_p \, \varphi'(0)^p \,  \sum_{k=N}^\infty (r_{k+1} -r_k)^p  \int_0^{2 \pi}\left| \frac{\varphi'(r_{k+1}
e^{i t})}{\varphi'(0)}\right|^p \, dt 
\end{eqnarray*}
for some constant $C_p$ depending only on $p$. If we assume for a moment that $
p> 2/5$ then for every positive integer $k$
$$\int_0^{2 \pi}\left| \frac{\varphi'(r_{k+1}
e^{i t})}{\varphi'(0)}\right|^p \, dt \le \frac{D_p}{(1-r_{k+1})^{3p-1}}\,
,$$
where  $D_p$ is some constant depending only on $p$, see
Theorem 1 in \cite{FM76}. 
Combining these results and using
$$ S_p:=C_p D_p \sup_{n } f_n'(0)^p $$
leads to 
\begin{eqnarray*}
\int \limits_{0}^{2 \pi} |\varphi(r e^{i t}) -\varphi(\rho e^{i t})|^p \,
dt &\le &
S_p   \sum_{k=N}^{\infty} \frac{(r_{k+1} -r_k)^p}{(1-r_{k+1})^{3p-1}} 
 \le  S_p \sum_{k=N}^{\infty} 
\frac{(r_{k+1} -r_k)^p}{(r-r_{k+1})^{3p-1}}\\[2mm] & \le &
S_p  \sum_{k=N}^{\infty} (r -r_{k+1})^{1-2p}
 \le  S_p\,   r^{1-2p}\,  \sum_{k=N}^{\infty}\left(
  \frac{1}{2^{k+1}} \right)^{1-2p}\\[2mm] &=& S_p\, \left(
  \frac{r}{2^{N+1}} \right)^{1-2p} \sum_{k=0}^{\infty}\left(
  \frac{1}{2^{1-2p}} \right)^{k} 
\\[2mm]
& \le & K_p\,  S_p \,  (r-r_{N+1})^{1-2p} 
\le K_p \,  S_p \,   (1-\rho)^{1-2p} \, ,
\end{eqnarray*}
where  $K_p$ is a constant depending only on $p$.
If $p\le 2/5$ choose $q >0$ such that $q\, p >2/5$ and apply H\"older's
inequality to arrive at the same estimate. In this way we obtain
\begin{equation*}
 \int \limits_{0}^{2 \pi} |f_n( e^{i t}) -f( e^{i t})|^p \,
dt \le 
\int \limits_{0}^{2 \pi} |f_n(\rho e^{i t}) -f(\rho e^{i t})|^p \,
dt + K (1-\rho)^{1-2p}
\end{equation*}
for some constant $K$. This finishes the proof of (i) $\Rightarrow$ (ii).
\hfill{$\blacksquare$}

\bigskip

{\bf Proof  of Lemma \ref{lem:B_Phi}  (c).}

Let $(g_n)_n \subset B_{\Phi}$ be a sequence which converges locally uniformly
in $\D$ to $g$ and let  $|g_n(z)|< K$ for all $n$. We note that $g\not \equiv
0$ since
$$\liminf_{|z| \to 1}
|g_n'(z)| \ge \inf_{|w|\le K} \Phi (w) >0 \, . $$  
So $g \in \mathcal{H}_0(\D)$ is univalent and bounded. 
Lemma \ref{lem:H^p} guarantees that we can subtract a
subsequence $(g_{n_j})_j$ such that $g_{n_j}(e^{i t}) \longrightarrow
g(e^{i t})$ a.\!\;e.\,. Thus the inequality
(\ref{eq:ex1}), which holds for every $g_{n_j} \in B_{\Phi}$, is valid also for
the limit function $g$ and Lemma \ref{lem:B_Phi} (b) implies
 $g \in B_{\Phi}$.~\hfill{$\blacksquare$}

\medskip

We next describe Beurling's geometric substitute of Perron's method for univalent supersolutions.
First, we define a ``reduced'' intersection of simply connected domains  
which is again a simply connected domain.

\begin{definition}[Reduced intersection]
Let $D_1$ and $D_2$ be two simply connected domains in $\C$ with $0 \in D_1
\cap D_2$. We denote by $\RI(D_1,D_2)$ the component of $D_1 \cap D_2$ which
contains the origin and call $\RI(D_1,D_2)$ the reduced intersection of $D_1$
and $D_2$. 
\end{definition}

Beurling's definition for the reduced intersection is formally different, but
equivalent. We note the following properties:
\begin{itemize}
\item[]{\bf (RI1)}  \, $\RI(D_1,D_2)$ is the largest simply connected domain which 
contains $0$ and is \phantom{\bf (RI1)\, \,\,}contained in $D_1 \cap D_2$.
\item[]{\bf (RI2)} \, If $D_1 \subseteq D_1'$ and $D_2 \subseteq D_2'$, then
$\RI(D_1,D_2) \subseteq \RI(D_1',D_2')$.
\item[]{\bf (RI3)} \, $\partial \RI(D_1,D_2) \subseteq \partial D_1 \cup
  \partial D_2$.
\end{itemize}

Now one can define an analogue of the Poisson modification for superharmonic functions.

\begin{definition}[Lower Beurling modification]
Let $g_1, g_2 \in \mathcal{H}_0(\D)$ be univalent functions. Then the  conformal
map $g$ from $\D$ onto the reduced intersection $\RI(g_1(\D),g_2(\D))$ normalized by $g(0)=0$ and
$g'(0)>0$ is called the lower Beurling modification of $g_1$ and $g_2$.
\end{definition}

\begin{lemma} \label{lem:lower}
Let $\Phi : \C \to \R$ be a positive, continuous and bounded function,
let $g_1, g_2$ be two univalent supersolutions for $\Phi$ and suppose
that $g_1$ is bounded. Then the lower Beurling
modification of $g_1$ and $g_2$ is also a univalent supersolution for $\Phi$.
\end{lemma}

Note that in contrast to the analogous statement for the upper Beurling
modification (Lemma \ref{lem:upper_mod}) we need to assume $g_1$ is bounded, but  we do not assume that $g_2$ is bounded.

\smallskip

{\bf Proof. } The proof is similar to the proof of Lemma \ref{lem:upper_mod}, so we only 
indicate what needs to be changed.
Let $\Psi : \C \to \R$ be a positive,
continuous and bounded function with $\Psi < \Phi$ and let 
$g$ denote the lower Beurling modification of $g_1$ and $g_2$, which maps
$\D$ conformally onto $D:=\RI(g_1(\D),g_2(\D))$. It suffices to show that $g
\in B_{\Psi}$.

\smallskip

 In view of Lemma \ref{lem:B_Phi} (d) and our assumption that $g_1$ is
bounded by some constant $L>0$, we see that for some $r_0 \in (0,1)$ the
functions $g_{1,r}(z):=g_1(
rz)$ belong to $B_{\Psi}$ for all $0<r_0<r<1$. Since $g_2$ may be unbounded we
cannot conclude that $g_{2,r}(z):=g_2(rz)$ also belong to
$B_{\Psi}$. However, for all $r<1$ sufficiently close to $1$ we have as a
partial substitute
\begin{equation*} \label{eq:ha}
|g_{2,r}'(\xi_0)| \ge \Psi(g_{2,r}(\xi_0)) \, \text{ for all } \, \xi_0 \in \partial \D
\, \text{ with } \, |g_{2,r}(\xi_0)| \le L \, . 
\end{equation*}
The proof of this  condition  is identical to the proof of Lemma
\ref{lem:extended0} (d).  We are now in a position to repeat part (i) of the  proof of Lemma \ref{lem:upper_mod} to show that the lower Beurling modification $g_r$ of $g_{1,r}$ and $g_{2,r}$ belongs to $B_{\Psi}$ for all $r$ sufficiently close to $1$. Finally, with obvious modifications of part (ii) of the proof
of Lemma \ref{lem:upper_mod}, we deduce $g_r \to g$, so $g \in B_{\Psi}$
by  Lemma \ref{lem:B_Phi} (c).
\hfill{$\blacksquare$}

\bigskip

{\bf Proof of Theorem  \ref{thm:mini_function}.}\\
Let $M:= ||\Phi||$ and $m:= \inf_{g \in B_{\Phi}}
g'(0)$. Choose a sequence $(g_n)_n \subset B_{\Phi}$ such that 
$g_n'(0)  \to m$. Since the function  $\hat{g}(z)=M z$ is a univalent
supersolution for $\Phi$ we can consider the lower Beurling modifications of $g_n$ and
$\hat{g}$, i.e.~the conformal maps
$\tilde{g}_n :\D \to {\rm RI}\, (g_n(\D), \hat{g}(\D))$
normalized by $\tilde{g}_n(0)=0$ and $\tilde{g}_n'(0)>0$. Note that  $\lim_{n \to \infty} \tilde{g}_n'(0)=m$ since
each $\tilde{g}_n$ is subordinate to $g_n$.
 By construction, the functions
$\tilde{g}_n$ form a bounded sequence of univalent supersolutions for $\Phi$. Thus there
is a subsequence $(\tilde{g}_{n_j})_j$ which converge locally uniformly in $\D$ to a
function $g^*$ with ${g^*}'(0)=m$. Further, $g^*$ belongs to $B_{\Phi}$,
by Lemma \ref{lem:B_Phi} (c), so $B^* \subset g^*(\D)$. 
To see that $ g^*(\D) \subset B^*$ 
we need to show $g^*(\D) \subseteq g(\D)$ for every $g \in B_{\Phi}$. Let 
$g \in B_{\Phi}$ and let $G$ be the lower Beurling modification of $g$ and $g^*$, so 
in particular $G(\D) \subseteq g^*(\D)$. On the other hand, $G$ belongs to 
$B_{\Phi}$, i.e., $G'(0) \ge m={g^*}'(0)$. Consequently, $G=g^*$.
Thus $g^*(\D)=G(\D) =\RI(g^*(\D),g(\D))\subseteq g(\D)$. \hfill{$\blacksquare$}

\section{Solutions} \label{sec:solutions}

\subsection{Existence of solutions for a generalized Beurling problem}
\label{par:gen}

In this paragraph we prove that there is always at least one 
solution $f \in \mathcal{H}_0(\D)$ for $\Phi$. In other words we will ensure
the existence of a
solution to the special Riemann--Hilbert--Poincar\'e problem
$$ \lim \limits_{z \to \xi} \left( |f'(z)|-\Phi(f(z)) \right)=0 \, , \qquad
\xi \in \partial \D \, .$$
 We wish to emphasize  
that this solution $f$ is not necessarily univalent, but we can find 
a locally univalent solution. The existence of at least one solution for
$\Phi$ (univalent or not)
will  be an indispensable ingredient for our proof below that the maximal univalent 
subsolution $f^* \in A_{\Phi}$ and also
the minimal univalent supersolution $g^* \in B_{\Phi}$ are actually  solutions for $\Phi$.
We note that Beurling does not use nonunivalent solutions for showing that the
  largest subsolution is in fact a solution. His reasoning, however, appears
  to be inconclusive to us.\footnote{This applies in particular to the properties of the
    auxiliary function $H$ stated on p.~126, lines 30 ff.~and p.~127, line 1
of \cite{Beu53}.} 

\medskip

In point of fact, one can even prescribe finitely many points $z_1, \ldots ,
z_n \in \D \backslash \{ 0\}$\footnote{Without loss of generality
we exclude $z=0$ as a possible
critical point, because of our normalization $f'(0)>0$ for any 
$f \in \mathcal{H}_0(\D)$} and always find a solution $f \in
\mathcal{H}_0(\D)$ with $z_1, \ldots,
z_n$ as its  critical points for a more general class of
Riemann--Hilbert--Poincar\'e problems. We begin with a characterization of the
solutions to such a generalized Beurling--type boundary value problem.

\begin{lemma} \label{lem:fix}
Let $\Phi : \partial \D \times \C \to \R$ be a positive, continuous and bounded
function and let  $f\in  {\cal H}_0(\D)$. 
Then $f$ is a solution for $\Phi$, i.e.
\begin{equation} \label{eq:rhp}
 \lim \limits_{z \to \xi} \left(|f'(z)|-\Phi\left(\xi,f(z)\right)\right)=0 \, , \qquad \xi \in \partial
 \D \, ,
  \end{equation}
if and only if $f$ extends continuously to $\overline{\D}$ and
\begin{equation} \label{eq:g}
f'(z)= B(z)
\exp \left( \frac{1}{2 \pi} \int
  \limits_{0}^{2 \pi} \frac{e^{it}+z}{e^{it}-z} \, \log \Phi\left(e^{it},f(e^{it})\right)
  \, dt \right)  \, , \quad  z \in \D\, ,
\end{equation}
where $B$ is a finite Blaschke product with $B(0)>0$.
\end{lemma}

{\bf Proof.}
Note that if $f$ is a solution for $\Phi$, then $f'$ is bounded on $\D$
by the maximum principle of subharmonic functions. Thus
$f$ is Lipschitz continuous on $\overline {\D}$ and by (\ref{eq:rhp}) $|f'|$ extends continuously
to $\overline{\D}$ (compare Lemma~\ref{lem:extended0} (a)).
 Since $|f'|$ is  positive on $\partial \D$, the function $f'$ has
only finitely many zeros in $\D\backslash \{0 \}$, say $z_1,
\ldots , z_n$. Set
$$B(z)= \eta \prod \limits_{j=1}^n \frac{z-z_j}{1-\overline{z_j} \, z} \, $$
with $\eta \in \partial \D$ such that $B(0)>0$. 
Then the function $z \mapsto \log(f'(z)/B(z))$ is holomorphic in $\D$ and 
$z \mapsto \log|f'(z)/B(z)|$ extends continuously to $\overline{\D}$ with
$ \log|f'(z)/B(z)|=\log \Phi(z, f(z))$ for $|z|=1$. Thus the Schwarz integral formula \cite[p.~42]{Pom92} leads to (\ref{eq:g}).

\smallskip

On the other hand, if $f$ extends continuously to $\overline{\D}$ and
$f'$ has the form (\ref{eq:g}), then  
it is clear that $f$
is a solution for $\Phi$.  \hfill{$\blacksquare$}

\begin{theorem} \label{thm:ma0}
Let $\Phi : \partial \D \times \C \to \R$ be a positive, continuous and bounded function and let
$z_1, \ldots, z_n \in \D \backslash \{ 0 \}$ be finitely many points.
Then there exists a solution $f \in \mathcal{H}_0(\D)$ for $\Phi$
with
 critical points $z_1, \ldots, z_n$ and no others.
\end{theorem}

\noindent
{\bf Proof.}
Let 
$$M:=\sup \limits_{(\xi,w) \in \partial \D \times \C} \Phi(\xi,w) \quad \text{ and } \quad
 B(z):=\eta \, \prod \limits_{j=1}^n \frac{z-z_j}{1-\overline{z_j} \, z} \, $$
with $\eta \in \partial \D$ such that $B(0)>0$. Further,
 let $\mathcal{H}(\D)$ denote the set of holomorphic functions $f : \D \to \C$.
Then
\begin{equation}\label{eq:setW}
 W:=\left\{ f \in \mathcal{H}(\D) \, \bigg| \, f(0)=0\, , \, f'(0) \ge 0 \, ,
\, |f'(z)| \le M \text{ for all } z \in \D \right\}
\end{equation} 
is a compact and convex subset of $\mathcal{H}(\D)$ (considered as a topological
vector space endowed with the topology of uniform convergence on compact
subsets of $\D$). Clearly, each $f \in W$ has a continuous extension to the
closure $\overline{\D}$, which we also denote by $f$. The operator ${\cal T} : W \to
{\cal H}(\D)$ given by
\begin{equation}\label{eq:operatorT}
 {\cal T}(f)(z):=\int \limits_{0}^{z} B(\zeta)
\exp \left( \frac{1}{2 \pi} \int
  \limits_{0}^{2 \pi} \frac{e^{it}+\zeta}{e^{it}-\zeta} \, \log \Phi(e^{it}, f(e^{it}))
  \, dt \right) \, d\zeta \, ,
\end{equation}
is then well--defined and it is not difficult to check that ${\cal T}$ is continuous. Note that
$ |{\cal T}(f)'(z)| \le M$ and
${\cal T}(f)'(0)>0$. Hence ${\cal T}(W) \subseteq W$. Thus, by Schauder's fixed point theorem,
  ${\cal T}$ has a fixed point $f_0 \in W \cap \mathcal{H}_0(\D)$ 
(which extends continuously to  $\overline{\D}$). By Lemma \ref{lem:fix}
this fixed point $f_0$ 
is a solution for $\Phi$ with critical points $z_1,
\ldots, z_n$ and no others. \hfill{$\blacksquare$}

\begin{corollary} \label{thm:ma}
Let $\Phi : \C \to \R$ be a positive, continuous and bounded function.
Then there exists a locally univalent solution $f \in \mathcal{H}_0(\D)$ for
$\Phi$.
\end{corollary}

\subsection{Regularity of solutions for the generalized
Beurling problem}\label{par:reg}

We now turn towards the boundary behaviour of the solutions $f\in {\cal
  H}_0(\D)$ to the generalized Beurling problem (\ref{eq:rhp}). For this we first introduce some notation. 

\medskip

For  $0<\alpha<1$ and $k=0,1,2, \ldots$ let $C^{k,\alpha}(\partial \D\times\C)$ denote
the set of all complex--valued functions on $\partial \D\times\C$ with all partial derivatives of
order $\le k$ continuous in $\partial \D\times\C$ and whose $k$--th order partial
derivatives are locally H\"older continuous with exponent $\alpha$ in
$\partial \D\times\C$. We say
a function $f \in \mathcal{H}_0(\D)$ belongs to the set $C^{k,\alpha}(\partial
 \D)$ if $f^{(k)}$ has a continuous extension to $\partial \D$ which is
H\"older continuous with exponent $\alpha$ on $\partial \D$ and to
the class $W_p^k(\partial \D)$, $0<p\le \infty$, if $f^{(k)}$ belongs to the
Hardy space $H^p$.

\begin{theorem}[Boundary regularity of solutions] \label{thm:regularity}
Let $\Phi : \partial \D\times\C \to \R$ be a positive, continuous and bounded
function and let $f$ be a solution to (\ref{eq:rhp}).
\begin{itemize}
\item[{\rm (a)}]
If $\Phi \in C^{k,\alpha}(\partial \D\times\C)$ for some $k\ge 0$ and $0<\alpha<1$,
then $f \in C^{k+1,\alpha}(\partial \D)$.

\item[{\rm (b)}]
If $\Phi \in C^{k}(\partial \D\times\C)$ for some $k\ge 1$,
then $f \in W^{k+1}_p(\partial \D)$ for all $0<p< \infty$.

\item[{\rm (c)}] 
 If $\Phi$ is real analytic on  $\{ (\xi, f(\xi))\, : \,
  \xi \in \partial \D \}$ 
then $f$ has an analytic extension across $\partial \D$. 
\end{itemize}
\end{theorem}

{\bf Proof.}
 We begin with some preliminary observations. If $\Phi$ belongs to
$C^{0,\alpha}(\partial \D \times \C)$ and $g$ is a Lipschitz continuous function on $\overline{\D}$,
then the function $\xi \mapsto \Phi(\xi, g(\xi))$ belongs to $C^{0,
  \alpha}(\partial \D)$. Let $k \ge 1$, then, if $\Phi \in
C^{k,\alpha}(\partial \D \times \C)$ and $g \in C^{k,
  \alpha}(\partial \D)$, the function $\xi \mapsto \Phi(\xi, g(\xi))$
 belongs to $C^{k,\alpha}(\partial \D)$.  
 The Herglotz integral of a function $g \in C(\partial \D)$ 
$$z \mapsto \frac{1}{2 \pi} \, \int 
  \limits_{0}^{2 \pi} \frac{e^{it}+z}{e^{it}-z} \, g(e^{it})\, dt $$
 belongs to $H^p$ for all $0<p<\infty$, see \cite[Chap.~II, Theorem
  3.1]{MG}. Moreover,  if $g \in C^1(\partial \D)$  standard Fourier techniques imply
\begin{equation}\label{eq:herglotz}
i z \left( \frac{1}{2 \pi} \, \int 
  \limits_{0}^{2 \pi} \frac{e^{it}+z}{e^{it}-z} \, g(e^{it})\, dt\right)'=
  \frac{1}{2 \pi} \, \int 
  \limits_{0}^{2 \pi} \frac{e^{it}+z}{e^{it}-z} \, \frac{d}{dt}\left(
  g(e^{it}) \right) \,  dt \, . 
\end{equation}

\smallskip

In the following let $f$ be a solution to (\ref{eq:rhp}). Then by Lemma \ref{lem:fix}, $f$ satisfies (\ref{eq:g}),
where $B$ is a finite Blaschke product with $B(0)>0$.
Recall, that $f$ has a Lipschitz
continuous extension to $\overline{\D}$ and $|f'|$ extends continuously to
$\partial \D$ (see the proof of  Lemma \ref{lem:fix}).  
 
\begin{itemize}
\item[(a)]
 If $\Phi \in C^{k,\alpha}(\partial \D \times \C)$, $k \ge 0$, then $\xi \mapsto
\log (\Phi(\xi,  f(\xi))$ belongs to
$C^{0,\alpha}(\partial \D)$, so formula (\ref{eq:g}) shows that
$f' \in C^{0,\alpha}(\partial \D)$ (see
\cite[Chap.~II, Corollary 3.3]{MG}). Hence $f \in C^{1,\alpha}(\partial \D)$.
Inductively, we obtain $f \in C^{k+1,\alpha}(\partial \D)$.

\item[(b)] If $\Phi \in C^k(\partial \D \times \C)$, $k \ge 1$, then  $\Phi
  \in C^{k-1,\alpha}(\partial \D \times \C)$ for every $0< \alpha <1$ and therefore
  $f$  belongs to $C^{k, \alpha}(\partial \D)$ by
  (a). So $\xi
  \mapsto\log  \Phi(\xi,
  f(\xi))$ belongs to $C^k(\partial \D)$, in particular to  $C^1(\partial
  \D)$. Combining (\ref{eq:herglotz}) and (\ref{eq:g}) yields
$$f''(z)=B'(z)\, \frac{f'(z)}{B(z)} +  f'(z) \frac{1}{iz} \left( \frac{1}{2 \pi} \, \int 
  \limits_{0}^{2 \pi} \frac{e^{it}+z}{e^{it}-z} \, \frac{d}{dt}\left(
  \log\Phi(e^{it}, f(e^{it}))\right) \,  dt   \right)\, ,\quad z \in \D\, .$$
Thus $f'' \in H^p$ for all $0<p<\infty$, which is the same as $f \in
  W^{2}_p(\partial \D)$ for $0<p<\infty$. Inductively, we obtain $f \in
  W_p^{k+1}(\partial \D)$, $0<p<\infty$.
\item[(c)]
If $\Phi$ does not explicitly depend on $z$, so $\Phi: \C \to \R$, then the
assertion is exactly Theorem 1.10 in \cite{Ro07}. The general case can be
handled in a similar way; we omit the details. \hfill{$\blacksquare$}
\end{itemize}

\begin{corollary}[Boundary regularity of univalent solutions]
  \label{thm:regularityunivalent}
Let $\Phi : \C \to \R$ be a positive, continuous and bounded function and let
$f \in \mathcal{H}_0(\D)$ be a univalent solution for $\Phi$. Then
$\Omega:=f(\D)$ is a Jordan domain bounded by a Lipschitz continuous Jordan curve.
\begin{itemize}
\item[(a)] If $\Phi \in C^{k,\alpha}(\partial \Omega)$,  
then $\partial \Omega$ is a Jordan curve of class
$C^{k+1,\alpha}$ for $0<\alpha<1$
and
$k=0,1, \ldots$.
\item[(b)] If $\Phi$ is real analytic on  $\partial \Omega$, then $\partial \Omega$ is 
  real analytic.
\end{itemize}
\end{corollary}

{\bf Proof.}
We just note that $f$ extends to a homeomorphism from $\overline{ \D}$ onto $\overline{f(\D)}$, see Lemma
\ref{lem:extended0} (a) and Lemma \ref{lem:B_Phi} (a).
\hfill{$\blacksquare$}

\medskip

The following example  shows that even though $|f'|$ has a continuous
extension to $\overline{\D}$ if $f$ is a  solution to Beurling's boundary value
problem,
this does not imply that $f'$ has a continuous extension to $\overline{\D}$.

\begin{example} \label{ex:pom}
Let $h$ be a conformal map from $\D$ onto the domain $D$ which is
obtained from  the rectangle $\{x+i y \, : \, -1 < x<1, \, -\pi/2<y<\pi/2\}$
by removing the vertical segements $1-1/k+i t$, $-\pi/2 < t \le 0$, $k=1,2, \ldots$.
Then  $\Re h$ is continuous up to $\partial \D$, while
$\Im h$ has no continuous extension to $\overline{\D}$, and $f \in
\mathcal{H}_0(\D)$ defined by $\log f'=h$ maps $\D$ onto a Jordan domain
$\Omega$ because $|\arg f'|=|\Im h|<\pi/2$, i.e.~$\Re f'>0$, 
cf.~{\cite[p.~45]{Pom92}}.
We define a continuous function $\Phi$ on $\partial \Omega$ by
$$ \Phi(w):=|f'(f^{-1}(w))|=\exp \left( \Re h(f^{-1}(w)) \right)$$
and extend $\Phi$ to a real--valued continuous bounded and non--vanishing
function on $\C$, so $f$ is a solution for $\Phi$, but $f'=e^h$ is certainly not continuous
up to $\partial \D$.
\end{example}

\subsection{The maximal univalent solution} \label{par:maxi}

\begin{theorem}\label{thm:uni_sol_A_Phi}
Let $\Phi:\C \to \R$ be a positive, continuous and bounded function. Then 
the maximal univalent subsolution $f^*$ for $\Phi$ is a solution to $\Phi$.
\end{theorem}
\noindent

{\bf Proof.}
Let $A^*=f^*(\D)$, let $U$ be the harmonic function in $A^*$ with boundary values $\log \Phi$
and let $\Psi(w):=\exp (U(w))$ inside $A^*$ and $\Psi(w):=\Phi(w)$ outside $A^*$.
Note that $f^* \in A_{\Psi}$.
We first show 
\begin{equation} \label{eq:help}
f^{* \, \prime}(0)=\max \limits_{h \in A_{\Psi}} h'(0) \, .
\end{equation}
By Theorem \ref{thm:maxi_function} there is a  maximal univalent function
$h^*$ 
in $A_{\Psi}$. In particular,
\begin{equation} \label{eq:help2}
 f^*(\D) \subseteq  h^*(\D) 
\end{equation}
and therefore $\partial h^*(\D) \subseteq \C \setminus A^*$.
Thus $\Psi=\Phi$ on $\partial h^*(\D)$, i.e.,  $\Psi(h^*(z))=\Phi(h^*(z))$ for $|z|=1$,
that is,
$$\limsup \limits_{|z| \to 1} \left( |h^{*\, \prime}(z)|-\Phi(h^*(z)) \right) \le 0 \, ,$$
so $h^* \in A_{\Phi}$. This implies $f^{*\, \prime}(0) \ge h^{*\, \prime}(0)$, which, combined with
(\ref{eq:help2}), gives $f^*=h^*$ and proves (\ref{eq:help}).

\smallskip

  By Corollary \ref{thm:ma} there is a
locally univalent holomorphic function $h_0 \in A_{\Psi}$ such that $|h_0'(z)|=\Psi(h_0(z))$ for
$|z|=1$. In view of Theorem \ref{thm:maxi_function}, $h_0(\D) \subseteq f^*(\D)$, so
$\log \Psi(h_0(z))$ is harmonic in $\D$ and coincides on $\partial \D$ with
the likewise harmonic function $\log |h_0'(z)|$. Thus
$|h_0'(z)|=\Psi(h_0(z))$ for every $z \in \D$. From (\ref{eq:help})
we therefore obtain
$$ f^{* \, \prime}(0)\ge h_0'(0)=\Psi(h_0(0))=\Psi(0)=\Psi(f^*(0)) \, .$$
We conclude that the harmonic function $\log  |f^{*\, \prime}(z)|-\log \Psi(f^*(z))
$, which is non--positive on $\partial \D$,
 must be  constant $0$. Since $\Phi=\Psi$ on $\partial A^*$, the proof of
 Theorem \ref{thm:uni_sol_A_Phi} is complete. \hfill{$\blacksquare$}

\subsection{The minimal univalent solution} \label{par:mini}

\begin{theorem}\label{thm:sol_B_Phi}
Let $\Phi:\C \to \R$ be a positive, continuous and bounded function. Then the
minimal univalent supersolution
$g^*$ for $\Phi$ is a solution to $\Phi$.

\end{theorem}

{\bf Proof.}
Let $g^*(\D)=B^*$. We define to each positive integer $n$ a  positive, continuous and bounded function
$\Phi_n: \C \to \R$ by
$$\Phi_n(w):= \frac{\Phi(w)}{n\, {\rm dist}(w, \,  \overline{B^*}) +1}\, , \quad
\text{where}  \quad {\rm dist}(w, \, \overline{B^*}):= \inf_{\eta \in \overline{B^*}}
|w-\eta|\, .$$
Note that $(\Phi_n)_n$ forms a monotonically decreasing sequence which is bounded above
by the function $\Phi$. Further, each $\Phi_n$ coincides with $\Phi$ on
$\overline{B^*}$  while $(\Phi_n)_n$ converges locally uniformly to $0$ on the
complement of $\overline{B^*}$.

\smallskip

In a first step we consider the sets $B_{\Phi_n}$. Let $g^*_n$ be the minimal  univalent
supersolution for $\Phi_n$ with $B_n^*:=g_n^*(\D)$,
see Theorem \ref{thm:mini_function}.
We will show that
$$g^*\equiv g_n^* $$ 
for all $n$.
Since $g^*(e^{it})
\in \partial B^*$ for a.e.~$e^{it} \in \partial \D$ and $g^*$ is a bounded univalent supersolution for $\Phi$ we
obtain by Lemma \ref{lem:B_Phi} (b)
$$
\log|g^{* \prime}(z)|  \ge  \frac{1}{2\pi} \int \limits_{0}^{2\pi} P(z, e^{i t})
\log \Phi(g^*(e^{i t}))\, dt = \frac{1}{2\pi} \int \limits_{0}^{2\pi} P(z, e^{i
  t})
\log \Phi_n(g^*(e^{i  t}))\, dt$$
for $z \in \D$ and all $n$. Thus $g^*$ is a univalent supersolution for every $\Phi_n$
which implies $B_n^* \subset B^*$ and $g_n^{* \, \prime}(0) \le g^{*\,
  \prime} (0)$.

\smallskip

On the other hand we see that every $g_n^*$ is also a univalent supersolution for $\Phi$ because 
$\partial B_n^* \subset\overline{ B^*}$ (see above) and so by Lemma
\ref{lem:B_Phi} (b)

$$
\log|g^{* \prime}_n(z)|  \ge  \frac{1}{2\pi} \int \limits_{0}^{2\pi} P(z, e^{i t})
\log \Phi_n(g_n^*(e^{i t}))\, dt = \frac{1}{2\pi} \int \limits_{0}^{2\pi} P(z, e^{i
  t})
\log \Phi(g^*_n(e^{i  t}))\, dt$$
for $z \in \D$ and all $n$.
This shows $B^* \subset B^*_n$ and $g^{* \, \prime}(0) \le g^{*\,
  \prime}_n(0)$. Combining these observations we deduce that  $g^*\equiv g^*_n$ for all $n$.

\smallskip
We now turn to  a study of the sets $A_{\Phi_n}$. Let $f_n^*$ be the maximal
univalent solution for $\Phi_n$ and $A_n^*:=f_n^*(\D)$,
see Theorem \ref{thm:maxi_function} and Theorem \ref{thm:uni_sol_A_Phi}.
The properties of $\Phi_n$ imply that
$A_{\Phi_{n+1}} \subseteq A_{\Phi_{n}} \subseteq \ldots \subseteq A_{\Phi}$
and therefore
$$
A^*_{n+1} \subseteq A_{n}^* \quad \text{ for all } n \ge 1 \, .
$$ 
Hence $(A_n^*)_n$ is a monotonically decreasing sequence of domains which converges to its kernel $A$,
i.e.~$A=\cap A_n^*$.  We note that 
$B^* \subset A$ since $B^* \subseteq A_n^*$ for all $n$ which in turn is a consequence of the
fact that $f_n^*$ is a univalent supersolution to $\Phi_n$. Now the corresponding conformal maps
$f_n^*$ from $\D$  onto $A_n^*$ converge locally uniformly to a conformal map
$f$ in ${\cal H}_0(\D)$ which maps $\D$ onto
$A$. Lemma \ref{lem:extended0} (c) tells us
that $f$ is a univalent subsolution to $\Phi$ as well as to $\Phi_n$ for every $n$.
In particular, $f$ is continuous on $\overline{\D}$.

\smallskip

Our next aim is to show that the boundary of $A$ is contained in the boundary of
$B^*$. For that it suffices to show that $\partial A =f(\partial \D) \subseteq
\overline{B^*}$. We assume to the contrary that there is a point $\zeta$ which
belongs to  $\partial A$ but not to $\overline{B^*}$. So there is a disk
$K(\zeta)$ about $\zeta$ such that 
$K(\zeta) \cap \overline{B^*}=\emptyset$.
The set $U:=K(\zeta ) \cap \partial A$ is open in $\partial A$ and if we
restrict the function $f$ to $\partial \D$ then the pre-image
$f^{-1}(U)$  is again an open set in $\partial \D$.  Pick $\tau \in f^{-1}(U)$ and
let $(z_j)_j \subseteq \D$ be an arbitrary sequence with $z_j \to \tau$. As $f$ is a
subsolution to every $\Phi_n$ we have
$$\limsup_{j \to \infty} |f'(z_j) | - \Phi_n(f(\tau)) \le 0\quad \text{for all }
n\, .$$ 
Since $f(\tau) \in \C \backslash \overline{B^*}$ and $\Phi_n(w) \to
0$ for $w \in \C \backslash \overline{B^*} $, we obtain by letting $n \to \infty$ 
$$\limsup_{j \to \infty}|f'(z_j) |=0\, .$$
In particular, the (angular) limit of $f'$ at $z=\tau$ exists and $=0$
for each $\tau \in U$, so $f' \equiv 0$ by Privalov's theorem
and the contradiction is apparent.

\smallskip

In the last step we will prove that $g^*\equiv f$. Taken this for granted 
$g^*$ is not only a univalent supersolution but also a subsolution to $\Phi$ and the
desired result follows. To observe that $g^*\equiv f$ we define a bounded,
positive and continuous function $\Psi: \C \to \R$ by
$$
\Psi(w)= 
\begin{cases}     
\exp(U(w)) \qquad &\text{for } w \in B^*\\[2mm]
\Phi(w)  \qquad &\text{for } w \in \C \backslash B^*\, ,
\end{cases}
$$
where $U$ is the harmonic function in $B^*$, which is  continuous on $\overline{B^*}$ 
 with boundary values $\log \Phi$.
Due to fact that $\partial A \subset \partial B^*$ we obtain by Lemma
\ref{lem:extended0} (b)
\begin{equation}\label{eq:sub1}
\log|f'(z)| \le 
\frac{1}{2\pi} \int_0^{2\pi} P(z, e^{it})\, \log\Phi(f(e^{it}))\,dt
=\frac{1}{2\pi} \int_0^{2\pi} P(z, e^{it})\,
\log\Psi(f(e^{it}))\,dt
\end{equation}
for $z \in \D$, in particular  $f$ is a univalent subsolution for $\Psi$.

\smallskip

Similarly, we deduce from Lemma \ref{lem:B_Phi} (b), since $g^*(e^{it}) \in \partial B^*$ for a.e.~$e^{it} \in \partial \D$,
\begin{equation}\label{eq:super1}
\log|g^{*\, \prime}(z)| \ge 
\frac{1}{2\pi} \int_0^{2\pi} P(z, e^{it})\, \log\Phi(g^*(e^{it}))\,dt
=\frac{1}{2\pi} \int_0^{2\pi} P(z, e^{it})\,
\log\Psi(g^*(e^{it}))\,dt
\end{equation}
for $z \in \D$, so $g^*$ is a  univalent supersolution for $\Psi$.
Furthermore, by inequality (\ref{eq:sub1}) and (\ref{eq:super1}) 

\begin{align*}
\log|f'(z)| \le U(f(z)) \quad \text{and} \quad \log|g^{*\, \prime}(z)| \ge U(g^*(z))
\end{align*}
for $z \in \D$.
Thus
$$\exp(U(0))=\exp(U(g^*(0))\le g^{* \, \prime}(0) \le f'(0) \le
\exp(U(f(0)))=\exp(U(0))\, ,$$
which in turn shows
$g^{* \, \prime}(0) = f'(0)$. Finally, by the principle of subordination, we
arrive at the conclusion $g^*\equiv f$.
\hfill{$\blacksquare$}

\section{Uniqueness and mapping properties of solutions} \label{sec:uniqueness}

Beurling \cite{Beu53} showed that the maximal and minimal univalent solution coincide
provided that $\log \Phi$ is superharmonic in $\C$, so there is only one
{\it univalent} solution $f \in \mathcal{H}_0(\D)$ in this case.
It is now easy to extend this uniqueness result to all {\it locally} univalent solutions.

\begin{theorem}
Let $\Phi : \C \to \R$ be a positive, continuous and bounded function and
suppose that $\log \Phi$ is superharmonic. Then
there is exactly one locally univalent solution for $\Phi$. This locally
univalent solution is univalent. 
\end{theorem}

{\bf Proof.}
Let $g$ be a locally univalent solution for $\Phi$ and
let $f^*$ denote the maximal univalent subsolution for $\Phi$, i.e.~$g(\D) \subseteq
f^*(\D)$, see  Theorem \ref{thm:maxi_function}. By  Theorem \ref{thm:uni_sol_A_Phi} $f^*$ is a univalent solution for $\Phi$. We need to show $g=f^*$.
The  functions
$u(z):=\log |{f^*}'(z)|$ and $v(z):=\log |g'(z)|$ are harmonic in $\D$ with 
$u(\xi)=\log \Phi(f^*(\xi))$ and $v(\xi)=\log \Phi(g(\xi))$ for all $\xi \in
\partial \D$. Since $\log \Phi(f^*(z))$ is superharmonic, we get
$u(z) \le \log \Phi(f^*(z))$ for all $z \in \D$. Now let
$w(z):={f^*}^{-1}(g(z))$, $z \in \D$, and consider the harmonic function $u
\circ w$ in $\D$. Then
$$ \limsup_{z \to \xi} u(w(z)) \le  \limsup_{z \to \xi} \log \Phi(f^*(w(z)))=
\lim_{z \to \xi} \log \Phi(g(z))=v(\xi) \, , \qquad \xi \in \partial \D \, .
$$
The maximum principle implies $u(w(0)) \le v(0)$, i.e., ${f^*}'(0) \le g'(0)$.
The uniqueness statement of Theorem \ref{thm:maxi_function} leads to the
conclusion that $g=f^*$. \hfill{$\blacksquare$}

\bigskip

The next theorem gives a condition on $\Phi$ which guarantess the uniqueness
of a solution with prescribed critical points for the extended Beurling
problem (\ref{eq:rhp}).

\begin{theorem}
Let $\Phi : \partial \D \times \C \to \R$ be a positive, continuous and
bounded function such that
\begin{equation}
|\Phi(\xi,w_1)-\Phi(\xi,w_2)|\le L\,|w_1-w_2|, \quad
\xi\in\partial \D,\ w_1,w_2\in\mathbb{C}
\end{equation}
for some positive constant $L$ with $L(1+M_0/m_0)<1$, where $M_0$ is chosen
such that
\begin{equation*}
\Phi(\xi,w)\le M_0 \, \, \,  \text{ for } \, |\xi|=1 \,,   |w|\le M_0 \quad
\text{ and } \quad 
m_0:=\min \{ \Phi(\xi,w):|\xi|=1,|w|\le M_0 \} \, .
\end{equation*}
Further, let $z_1, \ldots, z_n  \in \D \backslash \{ 0 \}$ be given. Then there exists a
unique solution $f \in {\cal H}_0(\D)$ for $\Phi$ with critical points $z_1,
\ldots, z_n$ and no others. 
\end{theorem}

{\bf Proof.}
We consider the set
 $W\subset {\cal H}(\D)$ defined by (\ref{eq:setW})  and the operator ${\cal T}:W \to {\cal
   H} (\D)$ given by (\ref{eq:operatorT}).
By Lemma \ref{lem:fix} a function  $f \in {\cal H}_0(\D)$ with critical points  $z_1,
\ldots, z_n$ is a solution for $\Phi$ if and only if $f$ is a fixed point of
the operator ${\cal T}$. 

\smallskip
  
In a first step we will show that $||f||\le M_0$ provided $f$ is a fixed point
of ${\cal T}$. For this we observe that
\begin{equation*}
\Phi(\xi,w) \le \Phi(\xi,M_0w/|w|) + L\,\big|w-M_0w/|w|\big| 
\le M_0 + L\,(|w|-M_0) < |w|
\end{equation*}
when $\xi \in \partial \D$ and $|w| >M_0$. Assume now that $f$ is a fixed point of
${\cal T}$ with $||f||>M_0$. Then 
\begin{equation*}
 ||{\cal T}f||\le \sup_{z \in \D} |B(z)| \, \exp \Bigg( \frac{1}{2 \pi} \int
  \limits_{0}^{2 \pi} P(z,e^{it}) \log \Phi(e^{it},f(e^{i t})) \,
  dt\Bigg)
\le\sup_{\xi \in \partial \D} \Phi(\xi, f(\xi)) <
  ||f|| 
\end{equation*}
which is a contradiction.

\smallskip

 For convenience we write
$$ {\cal S } f(z):= \frac{1}{2 \pi} \int
  \limits_{0}^{2 \pi} \frac{e^{it}+z}{e^{it}-z} \, \log \Phi(e^{it}, f(e^{it}))
  \, dt\, .$$
Now let $f_1$ and $f_2$ be two solutions for $\Phi$ both having critical
points $z_1, \ldots, z_n$. 
Then we obtain
\begin{equation*}
\begin{split}
&
 \left| \exp \left( {\cal S } f_2(z) \right)- \exp \left( {\cal S } f_1(z)\right)
 \right|\\[2mm]
&\le \left| \exp \left( \Re \left( {\cal S } f_2(z)\right)\right)- \exp\left(  \Re \left({\cal S } f_1(z)\right)\right)
 \right|+ M_0 \, \left| \exp \left( i \Im \left( {\cal S } f_2(z)\right)\right) - \exp
\left( i \Im \left({\cal S } f_1(z)\right) \right)
 \right|\\[2mm]
&\le \left| \exp \left( \Re \left( {\cal S } f_2(z)\right)\right) - \exp \left(\Re
 \left({\cal S } f_1(z)\right) \right)
 \right|+ M_0 \, \left|  \Im \left( {\cal S } f_2(z)\right) - \Im\left( {\cal S } f_1(z)\right) \right|\, .
\end{split}
\end{equation*}

Now applying the $|| \cdot||_2$ norm  of the Hardy
space $H^2$ to ${\cal T}f_2- {\cal T}f_1$ and using the facts
 that $||{\cal T} f_2-{\cal T} f_1||_2 \le ||\exp({\cal S} f_2)-\exp({\cal S}
f_1)||_2$ and
$||  \Im \left( {\cal S }
    f\right)||_2 \le ||\log \Phi(\cdot,f)||_2 $ (see \cite[p.~54]{Dur2000}) yields
\begin{equation*}
\begin{split}
||{\cal T}f_2 -{\cal T}f_1||_2 & \le ||\Phi(\cdot, f_2) -\Phi(\cdot, f_1)||_2 +M_0 ||\log
  \Phi(\cdot, f_2) -\log \Phi(\cdot, f_1)||_2  \\[2mm]
&\le L (1 + M_0 /m_0 )\,  ||f_2-f_1||_2 \, .
\end{split}
\end{equation*}
Since $L(1+M_0/m_0)<1$, we conclude
 $f_1\equiv f_2$.
\hfill{$\blacksquare$}

\bigskip
We now come to a uniqueness condition of different type.
It follows from  results of Gustafsson and Shahgholian  
\cite[Theorem 3.12]{GS} that if $\Phi : \C \to \R$ satisfies 
\begin{equation} \label{eq:tepper}
 \Phi(w) \le \frac{\Phi(\rho w)}{\rho} \, , \qquad 0<\rho<1 \, , \quad w \in
\C\, , 
\end{equation}
then every univalent solution $f \in \mathcal{H}_0(\D)$ for $\Phi$
maps $\D$ onto a starlike domain with respect to the origin.
The method of \cite{GS} is quite involved. It is based on a connection of Beurling's
boundary value problem with a class of free boundary value problem in PDE (see
Appendix 1 below for a discussion of this connection) and uses the ``moving
plane method''. Huntey, Moh and Tepper \cite{HMT03} showed that if strict inequality holds
in (\ref{eq:tepper}) for all $0<\rho<1$ and all $w \in \C$, then there is only
one univalent solution for $\Phi$. The following theorem gives more
precise information and includes the  results from \cite{GS} and
\cite{HMT03} mentioned above. It shows that under the condition (\ref{eq:tepper}) 
the maximal univalent solution $f^*$ for $\Phi$ is starlike and  
{\it every} univalent solution $f$ for $\Phi$  has the form 
$f=T f^*$ for some $T \in (0,1]$.

\begin{theorem} \label{thm:tep}
Let $\Phi : \C \to \R$ be a positive, continuous and bounded function which
satisfies (\ref{eq:tepper}). Then
the maximal univalent solution $f^*$ for $\Phi$ maps $\D$ onto  a starlike
domain with respect to 0 and there is a closed interval $I
\subset (0,1]$ such that the set of all univalent solutions for $\Phi$
is $\{T f^* \, : \, T \in I\}$. $f^*$ is the only univalent solution if in addition
\begin{equation} \label{eq:tepper2}
 \Phi(w)<\frac{\Phi(\rho w)}{\rho} \, , \qquad 0<\rho<1 \, , \quad w \in
\partial f^*(\D) \, . 
\end{equation}
\end{theorem}

{\bf Proof.} Note that $f \in A_{\Phi}$ implies $ T f \in A_{\Phi}$ for all
$0<T<1$. In particular, $A^*=f^*(\D)$ is starlike with respect to $w=0$.
Now let $f \in {\cal H}_0(\D)$ be a univalent solution for $\Phi$, i.e.,
$f(\D) \subseteq f^*(\D)$. Then there
 is a largest $T \in (0, 1]$ such that $T A^* \subseteq f(\D)$, so
the function $w(z):={f}^{-1}(T f^*(z))$ maps $\D$ into $\D$.
Since $f$ and $f^*$ map onto Jordan domains (see Corollary \ref{thm:regularityunivalent}),
the function $w$ extends continuously to $\overline{\D}$ and
$\xi_0:=w(\xi) \in \partial \D$ for at least some $\xi \in \partial \D$.
Note that $w(0)=0$ and $w'(0)>0$. 
  Assume $w \not\equiv \text{id}$. Then
the Julia--Wolff Lemma 
implies that $w$ has an angular derivative $w'(\xi)$ at $z=\xi$, where
$|w'(\xi)| \in (1,+\infty]$. Since 
$$|w'(z)| \le T \,\,  \frac{\sup_{z \in \D} |(f^*)'(z)|}{\inf_{z \in \D}
  |f'(w(z))|} \le \frac{T  \, ||\Phi||}{\inf_{z \in \D} |f'(w(z))|}  \le \frac{T  \, ||\Phi||}{\inf_{w \in f(\D)} |\Phi(w)|}  =:L <\infty \,
, \quad z \in \D\, ,$$
we have
$$ \left|\frac{w(\xi)-w(z)}{\xi-z}\right|=\lim \limits_{\eps \to 0}
\left| \frac{1}{z-\xi} \int \limits_{z}^{\xi-\eps} w'(s) \, ds \right| \le L
\, , \qquad z \in \D \, , 
$$
so $|w'(\xi)|$ is {\it finite}. Thus, if $\angle \lim$ denotes  the angular limit,
we obtain by using the Julia--Wolff lemma again
%
\begin{eqnarray*}
 |f'(\xi_0)| &=&  \angle \lim \limits_{z \to \xi}
 |f'(w(z))| 
 < \angle \lim \limits_{z \to \xi} |f'(w(z))|\cdot |w'(\xi)|= 
\angle \lim \limits_{z \to \xi} \left( |f'(w(z))|\cdot
|w'(z)|\right) \\ &=& \angle \lim \limits_{z \to \xi} T \, | {f^*}'(z)| =
T \, \Phi(f^*(\xi)) \le \Phi(T f^*(\xi))=\Phi(f(\xi_0)) \, ,
\end{eqnarray*}
a contradiction. Thus $f=T f^*$ for some $0<T \le 1$. In particular, there is
some $T_0$, $0 <T_0 \le 1$, such that $g^*= T_0
f^*$ where $g^*$ denotes the minimal univalent solution to $\Phi$. 

\smallskip

We now  prove that $T f^*$ is a solution to $\Phi$
for every $T \in I:=[T_0,1]$. Note that  for $\xi \in \partial \D$
\begin{eqnarray*}
 |(T f^*)'(\xi)|& =& T \Phi(f^*(\xi)) \le \Phi(T f^*(\xi))
=\Phi\left( \frac{T}{T_0} \, g^*(\xi)  \right) \le
\frac{T}{T_0} \Phi(g^*(\xi))=\frac{T}{T_0} |{g^*}'(\xi)| \\ & = & |(T f^*)'(\xi)| \,  . 
\end{eqnarray*}
Thus $T f^*$ is a solution for $\Phi$.

\smallskip

If (\ref{eq:tepper2}) holds and $f$ is a univalent solution for $\Phi$ different
from $f^*$, then $f=T f^*$ for some $T \in (0,1)$. But then
$$|f'(\xi)|= T \, |{f^*}'(\xi)|= T \Phi(f^*(\xi))<\Phi(T
f^*(\xi))=\Phi(f(\xi)) \, , \qquad \xi \in \partial \D \, ,$$
a contradiction.
\hfill{$\blacksquare$}

\smallskip

The following example illustrates the phenomena of
Theorem \ref{thm:tep}.

\begin{example}\label{ex:1}
Let 
$$ \Phi(w):=\begin{cases}  \sqrt{2|w|^2+1} \, & \, |w| \le 2 \, ,  \\
                           \quad 3 \, & \, 2<|w| \le 3 \, , \\
                           \quad |w| \, & \, 3<|w| \le 6 \, , \\
                           \quad 6   \, &  \, |w|>6 \, . 
\end{cases}$$
Then $\Phi$ is a positive, continuous and bounded function on $\C$ which satisfies condition
(\ref{eq:tepper}). Since $||\Phi||=6$  the maximal univalent solution is $f^*(z)=6z$.
Theorem \ref{thm:tep} implies that every univalent solution $f$ for $\Phi$ has the
form $f_r(z)=r \cdot z$ for some $0<r \le 6$. Now a direct calculation shows that 
$f_r(z)=r \cdot z$ is a univalent solution for $\Phi$ if and only if  $3 \le r
\le 6$. Thus the minimal univalent solution for $\Phi$ is $g^*(z)=3z$. We wish
to point out that $f(z)=z^2+z$ is a nonunivalent solution for $\Phi$ with $f(\D)
\subseteq g^*(\D)$.
\end{example}

\begin{remark}\label{rem:univalent_supersol}
The above example shows that the image domain of a holomorphic solution
for $\Phi$ is not necessarily contained in the image domain
of the corresponding minimal univalent solution. This is one of the reasons, 
why the class $B_{\Phi}$ is restricted to  univalent supersolutions.
In contrast, the image domain of every
solution is always contained in the image domain of the maximal univalent
solution. 

\end{remark}

\section{Appendix 1} \label{sec:appendix1}

We briefly discuss the relation of {\it univalent} solutions for Beurling's
boundary value problem with a class of free boundary value problems arising
in PDEs. 

\smallskip

Let $f \in \mathcal{H}_0(\D)$ be a {\it univalent} solution for $\Phi$.
Then $f$ maps $\D$ onto a simply connected domain $\Omega$ and
$u(w):=\log |f^{-1}(w)|$ is harmonic in $\Omega \backslash \{ 0\}$.
In fact, a quick computation shows that the pair $(u,\Omega)$ is a solution to the free boundary problem
\begin{equation} \label{eq:free}
\begin{array}{llll}
\Delta u(w) &= &2 \pi \delta_0(w)  &  w \in \Omega\\[3mm]
u(w) &=&0  \qquad \qquad \text{ for all } & w \in \partial \Omega\\[1mm]
|\nabla u(w)| &=& \displaystyle \frac{1}{\Phi(w)} &  w \in \partial \Omega \, .
\end{array}
\end{equation}
Here, $\delta_0$ denotes the Dirac delta function at $w=0$. Note that for
nonunivalent solutions $f$ for $\Phi$ the passage to (\ref{eq:free}) is not possible.
We call any pair $(u,\Omega)$, where $\Omega$ is a bounded domain in $\C$ and
$u : \Omega \to \R$  satisfies (\ref{eq:free}) a solution to (\ref{eq:free}).
It is easy to see that, if $(u,\Omega)$ is a solution to (\ref{eq:free}) with $\Omega$
simply connected, then there is an analytic function $F : \Omega \to \D$ with
$u(w)=\log |F(w)|$ with $F(0)=0<F'(0)$ and this function is actually a conformal map 
from $\Omega$ onto $\D$, so  $f:=F^{-1} \in \mathcal{H}_0(\D)$ is a solution
of Beurling's boundary problem (\ref{eq:boundarycondition}).

\smallskip

In particular, the regularity results in \cite{AC,KN77,GS} apply
immediately to any univalent solution $f \in {\cal H}_0(\D)$ for $\Phi$ and
show e.g.~that $f$ is in $C^{1,\beta}(\partial \D)$ for {\it some} $\beta \in
(0,1)$ whenever $\Phi \in C^{\alpha}(\C)$ for some $\alpha \in (0,1)$
and $f$ is real analytic, whenever $\Phi$ is real analytic. 
Theorem \ref{thm:regularity} is more precise and  more general as it also deals with
nonunivalent solutions for $\Phi$.

\smallskip

It was shown in \cite{GS,H} using PDE methods that there is
always a ``weak'' solution $(u,\Omega)$ to (\ref{eq:free}), i.e., $\Omega$ is
a bounded domain in $\C$, 
$u$ belongs
to the Sobolev space $H^1_0(\Omega)$ and the third condition in
(\ref{eq:free}) has to be interpreted in an appropriate weak sense, 
see \cite{H}. This information, however, does not suffice to guarantee that
there are  solutions
$f \in {\cal H}_0(\D)$ for $\Phi$, because the results in \cite{GS,H} do not
imply that there is a solution $(u,\Omega)$ to (\ref{eq:free}) where $\Omega$
is {\it simply connected}. In particular, Theorem \ref{thm:uni_sol_A_Phi} as well as
Theorem \ref{thm:sol_B_Phi} do not follow from the results in \cite{GS,H}.

\smallskip

\section{Appendix 2} \label{sec:appendix2}

In \cite[p.~120--121]{Beu53} Beurling studies sequences of bounded simply
connected domains $D \subset \C$ with $0 \in D$ using the standard
concept of kernel convergence. He prefers to speak of weak convergence instead
of kernel convergence. In point of fact, Beurling uses Pommerenke's definition
\cite[p.~13]{Pom92} of kernel convergence which is only formally different 
from the usual definition.
An important role in Beurling's approach is played by strictly shrinking
sequences of simply connected domains $D_n$, i.e., $\overline{D_{n+1}} \subset
D_n$ for all $n=1,2, \ldots$.
Beurling calls a domain $D \subseteq \C$ a domain of Schoenfliess type, if
 $\Omega:=\hat{\C} \backslash
\overline{D}$ is simply connected and $\partial D \subseteq \partial \Omega$.
He asserts (\cite[p.~121]{Beu53}) that if a  strictly shrinking sequence
of simply connected domains $D_n$ with $w_0 \in D_n$ 
converges (weakly or in the sense of kernel
convergence with respect to $w_0$) to a domain $D$ with $w_0 \in D$, then
$D$ is necessarily of Schoenfliess type. This assertion is repeatedly used by
Beurling and turns out to be 
particularly important in his proof that the minimal univalent supersolution
is a solution (see \cite[p.127--130]{Beu53}). However, the kernel of a 
strictly shrinking sequence of simply connected domains is not necessarily
of Schoenfliess type as the following example shows.

\begin{example} \label{ex:gegen}
Let $D$ be the bounded simply connected domain as shown on the left side of
Figure \ref{pic:1}. Note that  $\hat{\C} \backslash \overline{D}$ consists
of two components, the inner disk and an unbounded component. Thus $D$ is not
of Schoenfliess type. However, $D$ can be obtained as the kernel (with respect
to any point $w_0 \in D$) of simply connected domains $D_n$, where
$\overline{D_{n+1}} \subset D_n$. The domains $D_1$ (blue), $D_2$ (red) and
$D_3$ (green) are displayed on the right side of Figure \ref{pic:1}.
\vspace{-0.9cm}
\begin{figure}[h]
\hspace*{1cm}\epsfig{file=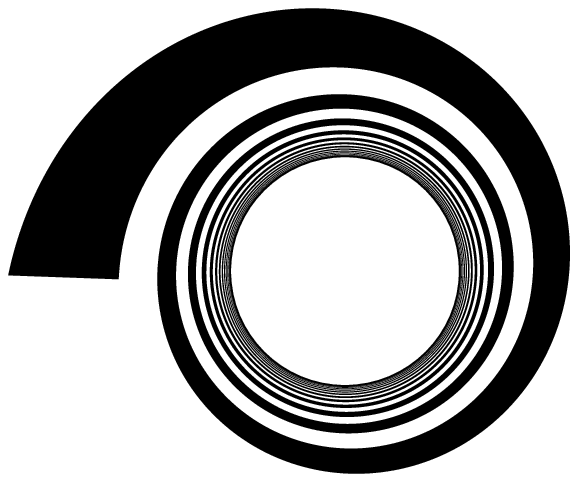,width=7cm,angle=90} \hspace*{2cm} \epsfig{file=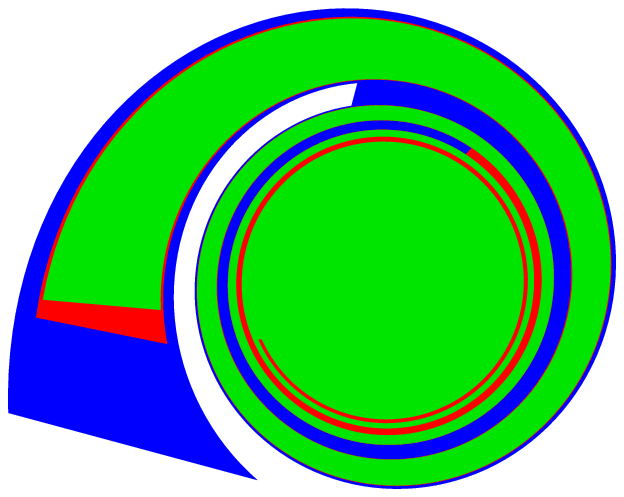,width=7cm,angle=90}
\vspace*{-1.5cm}
\caption{Strictly shrinking domains which converge to a domain
not of Schoenfliess type} \label{pic:1}
\end{figure}
\end{example}

\noindent
{\it Address:}

\noindent
\begin{minipage}[t]{8cm}
Florian Bauer, Daniela Kraus and Oliver Roth\\
Institut f\"ur Mathematik\\
Universit\"at W\"urzburg\\
Am Hubland \\
97074 W\"urzburg\\
Germany
\end{minipage} \hspace{1.7cm}
\begin{minipage}[t]{8cm}
Elias Wegert, \\
Institut f\"ur Angewandte Analysis\\  
TU Bergakademie Freiberg\\
Pr\"uferstr.\,9\\ 09596 Freiberg \\
Germany
\end{minipage}

\end{document}